\documentclass[11pt]{amsart}

\usepackage{amsmath}
\usepackage{amssymb}
\usepackage{amsthm}
\usepackage{mathtools}
\usepackage{bbm}
\usepackage[normalem]{ulem}
\usepackage[dvipsnames]{xcolor}
\usepackage[
  colorlinks=true,
  linkcolor=blue,
  citecolor=blue,
  urlcolor=blue,
  pdftitle={Submatrices of character tables and centers of q-Schur superalgebras},
  pdfauthor={Li Luo and Senlei Yu}
]{hyperref}

\allowdisplaybreaks
\numberwithin{equation}{section}

\theoremstyle{plain}
\newtheorem{thm}{Theorem}[section]
\newtheorem{lem}[thm]{Lemma}
\newtheorem{prop}[thm]{Proposition}

\theoremstyle{definition}

\theoremstyle{remark}
\newtheorem{rmk}[thm]{Remark}

\DeclareMathOperator{\End}{End}

\DeclareMathOperator{\Span}{span}

\DeclareMathOperator{\rank}{rank}

\newcommand{\id}{\operatorname{id}}

\newcommand{\nc}{\newcommand}
\nc{\redtext}[1]{\textcolor{red}{#1}}
\nc{\red}[1]{\redtext{ #1}}
\nc{\orange}[1]{\textcolor{orange}{#1}}
\nc{\browntext}[1]{\textcolor{brown}{#1}}
\nc{\greentext}[1]{\textcolor{green}{#1}}
\nc{\bluetext}[1]{\textcolor{blue}{#1}}
\nc{\brown}[1]{\browntext{ #1}}
\nc{\green}[1]{\greentext{ #1}}
\nc{\blue}[1]{\bluetext{ #1}}

\title[Centers of $q$-Schur superalgebras]
{Submatrices of character tables and centers of ordinary and queer $q$-Schur superalgebras}

\author[Li Luo]{Li Luo}
\author[Senlei Yu]{Senlei Yu}
\address{School of Mathematical Sciences, Key Laboratory of MEA (Ministry of Education) \& Shanghai Key Laboratory of PMMP, East China Normal University, Shanghai 200241, China}
\email{lluo@math.ecnu.edu.cn (Luo)}
\email{52280155006@stu.ecnu.edu.cn (Yu)}

\date{}

\begin{document}

\begin{abstract}
We obtain bases of the centers of the ordinary $q$-Schur superalgebra $\mathcal{S}_q(m|n,r)$ and the queer $q$-Schur superalgebra $\mathcal{Q}_q(n,r)$ in terms of Casimir elements of the Hecke algebra $\mathcal{H}_r$ and the Hecke-Clifford algebra $\mathcal{H}_r^c$, respectively. The key step is to find an invertible submatrix of the character table of $\mathcal{H}_r$ or $\mathcal{H}_r^c$ whose column index set is the support of the corresponding $q$-tensor module appearing in Schur-Weyl duality: for $\mathcal{H}_r$ this is the set $H_{m,n}(r)$ of $(m|n)$-hook partitions of $r$, while for $\mathcal{H}_r^c$ this is the set $\mathrm{SP}_n(r)$ of strict partitions of $r$ of length $\leq n$. Unlike the case of ordinary $q$-Schur algebras, the submatrix formed by taking the row index set to be the same as the column index set need not be invertible. We determine that the row index set should be chosen as $H_{\min\{m,n\},\max\{m,n\}}(r)$ for $\mathcal{S}_q(m|n,r)$, and as the preimage of $\mathrm{SP}_n(r)$ under the Glaisher correspondence for $\mathcal{Q}_q(n,r)$. 
\end{abstract}

\maketitle


\section{Introduction}

An algebra and the endomorphism algebra of a module over it are always closely related. The most typical example is the Schur-Weyl duality between the symmetric group (resp. its $q$-deformation, the Hecke algebra) and the Schur algebra (resp. the $q$-Schur algebra). This duality closely links the representation theories of the two algebras. Moreover, the structures of the two algebras themselves are also highly interrelated. For instance, in Fu's paper \cite{Fu21}, the center of the $q$-Schur algebra is formulated in terms of the center of the Hecke algebra.

This principle for centers of algebras is also valid for centers of superalgebras. Let $A=A_{\bar{0}}\oplus A_{\bar{1}}$ be a finite-dimensional split semisimple superalgebra over a field $F$ of $\mathrm{char} F\neq2$, and $M$ a nonzero finite-dimensional right $A$-supermodule. Let $\{e_\lambda\mid \lambda\in\Lambda\}$ be the primitive even central idempotents of $A$, where $\Lambda$ is an index set. Put $\Lambda_{\bar{1}}=\{\lambda\in\Lambda\mid Ae_\lambda \text{ is of type $\mathsf{Q}$}\}$. For each $\lambda\in\Lambda_{\bar{1}}$, choose $\omega_\lambda\in Z(Ae_\lambda)_{\bar{1}}$ with $\omega_\lambda^2=e_\lambda$. Then $\{e_\lambda\mid\lambda\in\Lambda\}$ and $\{\omega_\lambda\mid\lambda\in\Lambda_{\bar{1}}\}$ form bases of $Z(A)_{\bar{0}}$ and $Z(A)_{\bar{1}}$, respectively. Given a finite-dimensional right $A$-module $M$, denote by $\Lambda_M=\{\lambda\in\Lambda\mid Me_\lambda\neq0\}$ the support of $M$. Let $B=\mathrm{End}_A(M)$ be the endomorphism superalgebra of $M$ over $A$. We show the following relationship between $Z(A)$ and $Z(B)$ in Theorem~\ref{prop:supp-center}. 
The homomorphism $\eta_M: A^{\mathrm{op}} \to \mathrm{End}_F(M)$ induces a surjective graded algebra homomorphism $Z(A)\to Z(B)$ with $Z(B)_{\bar{0}}=\bigoplus_{\lambda\in\Lambda_M} F \eta_M(e_\lambda)$ and $Z(B)_{\bar{1}}=\bigoplus_{\lambda\in\Lambda_M\cap \Lambda_{\bar{1}}} F \eta_M(\omega_\lambda)$.
Consequently, if $(\xi_\mu)_{\mu\in I}$ is a family in $Z(A)_{\bar{0}}$ with $\xi_\mu=\sum_{\lambda\in\Lambda} c_{\mu\lambda}e_\lambda$, then $(\eta_M(\xi_\mu))_{\mu\in I}$ is a basis of $Z(B)_{\bar{0}}$ if and only if the restricted coefficient matrix $(c_{\mu\lambda})_{\mu\in I,\lambda\in\Lambda_M}$ is an invertible square matrix. 

Let $A$ be the Hecke algebra $\mathcal{H}_r$ or the Hecke-Clifford algebra $\mathcal{H}_r^c$ over a suitable field $\mathbb{K}$ such that $A$ is split semisimple. We regard $\mathcal{H}_r$ as a purely even superalgebra, and $\mathcal{H}_r^c$ as a superalgebra whose Clifford generators are odd. They admit Schur-Weyl dualities with the quantum general linear supergroup $\mathbf{U}_q(\mathfrak{gl}_{m|n})$ (see \cite{Mi06}) and the quantum queer supergroup $\mathbf{U}_q(\mathfrak{q}_n)$ (see \cite{Ol92}) on the $q$-tensor superspaces $(\mathbb{K}^{m|n})^{\otimes r}$ and $(\mathbb{K}^{n|n})^{\otimes r}$, respectively. The resulting endomorphism superalgebras in these dualities are the ordinary $q$-Schur superalgebra $\mathcal{S}_q(m|n,r)$ and the queer $q$-Schur superalgebra $\mathcal{Q}_q(n,r)$. If we have a basis $(\xi_\mu)$ of $Z(\mathcal{H}_r)_{\bar{0}}=Z(\mathcal{H}_r)$ (resp. $Z(\mathcal{H}_r^c)_{\bar{0}}$) which is a $q$-deformation of the class-sum basis of $Z(\mathfrak{S}_r)$ (resp. $Z(\mathfrak{S}^c_r)$), then the coordinate matrix of $(\xi_\mu)$, with respect to the basis consisting of the primitive central idempotents, is a $q$-deformation of the character table by the corresponding Frobenius formula. Therefore, we can obtain a basis of $Z(\mathcal{S}_q(m|n,r))_{\bar{0}}=Z(\mathcal{S}_q(m|n,r))$ (resp. $Z(\mathcal{Q}_q(n,r))_{\bar{0}}$) in terms of $(\xi_\mu)$ after finding an invertible submatrix of the character table corresponding to the support of the associated $q$-tensor superspace.

In the non-super case treated in \cite{Fu21}, Fu takes $(\xi_\mu)_{\mu\in\mathrm{Par}(r)}$ to be the Geck-Rouquier basis (cf. \cite{GR97}) or the Jones basis (cf. \cite{Jo90}) of $Z(\mathcal{H}_r)$ to capture a basis of the center of the $q$-Schur algebra $\mathcal{S}_q(n,r)=\mathrm{End}_{\mathcal{H}_r}((\mathbb{K}^n)^{\otimes r})$. Here the two bases are both indexed by the set $\mathrm{Par}(r)$ of partitions of $r$, as well as the basis $\{e_\lambda\mid\lambda\in\mathrm{Par}(r)\}$ consisting of the primitive central idempotents. The support of $(\mathbb{K}^n)^{\otimes r}$ for $\mathcal{H}_r$ is the set $\mathrm{Par}_n(r)$ of partitions of $r$ with at most $n$ parts. Fu shows that the restricted submatrix $(c_{\mu\lambda})_{\mu,\lambda\in\mathrm{Par}_n(r)}$ of the character table is invertible and therefore obtains a basis of $Z(\mathcal{S}_q(n,r))$ in terms of $\xi_\mu$, ($\mu\in \mathrm{Par}_n(r)$). Note that for the restricted submatrix $(c_{\mu\lambda})_{\mu,\lambda\in\mathrm{Par}_n(r)}$, the index sets of rows and of columns coincide.

In the super case, we know that the associated support for $\mathcal{H}_r$ (resp. $\mathcal{H}_r^c$) is the set $H_{m,n}(r)$ of $(m|n)$-hook partitions of $r$ (resp. the set $\mathrm{SP}_n(r)$ of strict partitions of $r$ with at most $n$ parts). 
However, for $\mathcal{H}_r$, the associated submatrix of its character table, with rows and columns both indexed by the supports, may not be invertible; while for $\mathcal{H}_r^c$, the index sets of rows and of columns of its character table are different to begin with.
It is a key point of this paper that, for $\mathcal{S}_q(m|n,r)$, we should take the associated restricted submatrix whose row and column index sets are $H_{\min\{m,n\},\max\{m,n\}}(r)$ and $H_{m,n}(r)$, respectively (cf. Lemma~\ref{thm:gl-mat-ns}); while for $\mathcal{Q}_q(n,r)$, we should take the one whose row and column index sets are $\operatorname{OP}_n^G(r)$ and $\mathrm{SP}_n(r)$, respectively (cf. Lemma~\ref{thm:q-mat-class}). Here $\operatorname{OP}_n^G(r)$ is the set of those odd partitions corresponding to $\mathrm{SP}_n(r)$ via the Glaisher correspondence. By these facts, we obtain a basis of $Z(\mathcal{S}_q(m|n,r))_{\bar{0}}=Z(\mathcal{S}_q(m|n,r))$ (resp. $Z(\mathcal{Q}_q(n,r))_{\bar{0}}$) in terms of the Casimir elements of $Z(\mathcal{H}_r)$ (resp. $Z(\mathcal{H}_r^c)$) in Theorem~\ref{thm:gl-center-bases} (resp. Theorem~\ref{thm:q-even-bases}). It also makes sense to get $Z(\mathcal{S}_q(m|n,r))$ through the Geck-Rouquier basis and the Jones basis of $Z(\mathcal{H}_r)$, which was recently done by Fu-Luo-Sun \cite{FLS26}. Furthermore, we also give an explicit basis of $Z(\mathcal{Q}_q(n,r))_{\bar{1}}$ in Theorem~\ref{thm:q-even-bases} to complete a basis of $Z(\mathcal{Q}_q(n,r))$.

The paper is organized as follows. In Section~\ref{sec:comb}, we recall some combinatorics about partitions and symmetric functions. In Section~\ref{sec:quantum}, we recall several algebras and superalgebras, including Hecke algebras, Hecke-Clifford algebras, ordinary $q$-Schur superalgebras, and queer $q$-Schur superalgebras. Section~\ref{sec:center} is the main text of this paper, in which we set up a general relationship between the centers of a split semisimple superalgebra $A$ and an endomorphism superalgebra $\mathrm{End}_A(M)$ of an $A$-supermodule $M$, and then apply it to get centers of ordinary $q$-Schur superalgebras, and queer $q$-Schur superalgebras through the centers of Hecke algebras and Hecke-Clifford algebras, respectively.

\subsection*{Convention}
Throughout this paper, all non-quantized vector spaces and algebras are defined
over the complex field $\mathbb C$. Let $v$ be an indeterminate and set $q=v^2$. Quantized spaces and algebras are defined over the algebraically closed field $\mathbb{K}=\overline{\mathbb{C}(v)}$. Here we take the ground field $\overline{\mathbb{C}(v)}$ to ensure split semisimplicity of the Hecke-Clifford algebras and the queer $q$-Schur superalgebras, although it suffices to work over the field $\mathbb{C}(q)$ when we consider the Hecke algebras and the ordinary $q$-Schur algebras.

For every algebra $A$, let $Z(A)$ denote the center of its underlying
associative algebra. If $A=A_{\bar0}\oplus A_{\bar1}$ is a
superalgebra, let $Z(A)_{\bar i}=Z(A)\cap A_{\bar i}$ for
$\bar i\in\mathbb Z_2$. Thus $Z(A)$ always denotes the ordinary
center, not the supercenter.

\subsection*{Acknowledgments}
The work is partially supported by the National Key R\&D Program of China (No. 2024YFA1013802) and the NSF of China (Nos. 12371028 and 12631002).

\section{Combinatorics} \label{sec:comb}

\subsection{Partitions}
Let $r$ be a positive integer. A partition $\lambda\vdash r$ of $r$ is a weakly decreasing sequence
$\lambda=(\lambda_1,\lambda_2,\ldots)$ of nonnegative integers with
$\sum_{i\geq1}\lambda_i=r$.
We denote its size by $|\lambda|=r$ and its length by
$\ell(\lambda)=\#\{i~|~\lambda_i>0\}$. Hence, $\lambda_i=0$ for $i>\ell(\lambda)$. Denote by $\operatorname{Par}(r)$ the set of partitions of $r$. We write
$\lambda'$ for the conjugate partition of $\lambda$. 

For any two partitions $\mu\vdash r_1$ and $\delta\vdash r_2$, let $\mu\cup\delta\vdash r_1+r_2$
denote the partition obtained by rearranging all their positive
parts, with multiplicities retained, in weakly decreasing order.

\subsubsection{Hook partitions}
Let $m,n$ be nonnegative integers with $m+n>0$. A partition $\lambda\vdash r$ is called an $(m|n)$-hook partition if
$\lambda_{m+1}\leq n$. Denote
$$
H_{m,n}(r)=\{\lambda\vdash r:\lambda_{m+1}\leq n\},
\quad
B_{m,n}(r)=H_{\min\{m,n\},\max\{m,n\}}(r).
$$

\begin{lem}
For every $\lambda\vdash r$, one has
$\lambda\in H_{m,n}(r)$ if and only if
$\lambda'\in H_{n,m}(r)$. Consequently,
$|H_{m,n}(r)|=|H_{n,m}(r)|$ and
\begin{equation}\label{eq:BH-card}
|B_{m,n}(r)|=|H_{m,n}(r)|.
\end{equation}
\end{lem}

\begin{proof}
The condition $\lambda_{m+1}\leq n$ says that the Young diagram of
$\lambda$ does not contain the box $(m+1,n+1)$. After transposition,
this is equivalent to $\lambda'_{n+1}\leq m$. If $m\leq n$, then
$B_{m,n}(r)=H_{m,n}(r)$. If $m>n$, conjugation gives a bijection
$H_{n,m}(r)\to H_{m,n}(r)$, which proves the cardinality statement.
\end{proof}

\subsubsection{Strict and odd partitions}
Let $\operatorname{SP}(r)$ be the set of strict partitions of $r$, i.e. the partitions whose positive parts are pairwise distinct.
Let $\operatorname{OP}(r)$ be the set of partitions of $r$ whose positive parts are all odd. 

There is a bijection $G:\operatorname{OP}(r)\to\operatorname{SP}(r)$,
called the Glaisher correspondence, as follows.
Take
$\mu=(1^{m_1}3^{m_3}5^{m_5}\cdots)\in\operatorname{OP}(r)$. Write each multiplicity in its binary expansion as $m_j=\sum_{a\geq0}\epsilon_{j,a}2^a$ with
$\epsilon_{j,a}\in\{0,1\}$. The partition $G(\mu)$ contains the
part $2^aj$ precisely when $\epsilon_{j,a}=1$. Conversely, a part
$2^{a}j$, with $j$ odd, is replaced by $2^a$ copies of $j$.

For $n\geq1$, denote
$$
\operatorname{SP}_n(r)
=\{\lambda\in\operatorname{SP}(r)~|~\ell(\lambda)\leq n\}
\quad \text{and} \quad
\operatorname{OP}_n^G(r)=G^{-1}(\operatorname{SP}_n(r)),
$$
so that
\begin{equation}\label{eq:q-index-card}
|\operatorname{OP}_n^G(r)|=|\operatorname{SP}_n(r)|.
\end{equation}

Finally, we denote
$$
\operatorname{SP}^{-}(r):=\{\lambda\in\operatorname{SP}(r)\mid
\ell(\lambda)\text{ is odd}\}\quad\text{and}\quad
\operatorname{SP}_n^{-}(r)
=\operatorname{SP}^{-}(r)\cap\operatorname{SP}_n(r).
$$

\subsubsection{Orders on partitions}
For a partition $\alpha$ and $k\in\mathbb{N}$, let $S_k(\alpha)$ denote the sum of the first $k$ parts of $\alpha$; it is understood that $S_0(\alpha)=0$. The  dominance order $\unrhd$ on $\operatorname{Par}(r)$ is defined by, for $\alpha,\beta\in\operatorname{Par}(r)$,
$$\alpha\unrhd\beta \quad\text{iff}\quad S_k(\alpha)\geq S_k(\beta), \quad \forall k\geq 1.
$$
We write $\alpha\triangleright\beta$ if $\alpha\unrhd\beta$ and
$\alpha\neq\beta$.

\begin{lem}
\label{lem:perm-major}
Let $\gamma=(\gamma_1,\ldots,\gamma_t)$ be a partition with
$\gamma_t\geq t$.  For a permutation $\sigma\in\mathfrak S_t$ of $\{1,2,\ldots,t\}$, let $\gamma_\sigma$ be the  weakly decreasing rearrangement
of the positive integers $\gamma_i-i+\sigma(i)$, ($i=1,2,\ldots,t$), which is a partition with size $|\gamma|$. Moreover,
$\gamma_\sigma\unrhd\gamma$, with equality if and only if
$\sigma=\mathrm{id}$.
\end{lem}

\begin{proof}
It is clear that $\gamma_\sigma$ is a partition with size $|\gamma|$ because $\sum_{i=1}^ti=\sum_{i=1}^t\sigma(i)$.

For $1\leq k\leq t$, we have 
$$
\begin{aligned}
S_k(\gamma_\sigma)\geq
\sum_{i=1}^k\bigl(\gamma_i-i+\sigma(i)\bigr)=
\sum_{i=1}^k\gamma_i
-\frac{k(k+1)}2+\sum_{i=1}^k\sigma(i)\geq
\sum_{i=1}^k\gamma_i=S_k(\gamma),
\end{aligned}
$$
which implies
$\gamma_\sigma\unrhd\gamma$.

If $\gamma_\sigma=\gamma$, then the equality must hold
throughout the preceding inequalities for every $k$, i.e.,
$\sum_{i=1}^k\sigma(i)=k(k+1)/2$ for every $k$. It follows inductively that $\sigma(k)=k$ for all $k$, i.e., $\sigma=\mathrm{id}$. Conversely, it is clear that $\gamma_{\sigma}=\gamma$ if $\sigma=\mathrm{id}$.
\end{proof}

\begin{lem}
\label{lem:union-dom}
If $\nu\triangleright\gamma$, then
$\nu\cup\delta\triangleright\gamma\cup\delta$ for every partition
$\delta$.
\end{lem}

\begin{proof}
By induction on $\ell(\delta)$, it suffices to deal with $\delta=(a)$ for an arbitrary positive
integer $a$.
Clearly,
$$
S_k(\rho\cup(a))
=\max\{S_k(\rho),\,S_{k-1}(\rho)+a\}
\qquad(\forall k\geq1).
$$
Since $\nu\unrhd\gamma$, this gives
$S_k(\nu\cup(a))\geq S_k(\gamma\cup(a))$ for every $k\geq1$.
Hence $\nu\cup(a)\unrhd\gamma\cup(a)$.
If $\nu\cup(a)=\gamma\cup(a)$, deleting one occurrence of $a$
from both sides gives $\nu=\gamma$, contrary to
$\nu\triangleright\gamma$.
\end{proof}

Moreover, we will use another partial order $\succ$ on $\operatorname{Par}(r)$ defined by
\begin{equation}\label{def:ps-order}
\alpha\succ\beta \quad\text{iff} \quad  \ell(\alpha)<\ell(\beta), \ \text{or},\ 
\ell(\alpha)=\ell(\beta)\ \text{and}\ 
\beta\triangleright\alpha.
\end{equation}

\subsection{Symmetric groups}
The symmetric group $\mathfrak S_r$, consisting of permutations of $\{1,2,\ldots,r\}$, is generated by simple reflections $s_1,s_2,\ldots, s_{r-1}$. For $\mu\vdash r$, let $C_\mu$ be the conjugacy class of cycle type $\mu$, and let
\begin{equation}\label{def:cmu}
c_\mu=\sum_{w\in C_\mu}w\in Z(\mathbb C\mathfrak S_r)
\end{equation}
be its class sum. Write $m_i(\mu)$ for the
multiplicity of $i$ in $\mu$, and denote
$$z_\mu=\prod_{i\geq1}i^{m_i(\mu)}m_i(\mu)!,$$ then
$|C_\mu|=r!/z_\mu$. For $\lambda\vdash r$, let $S^\lambda$ be
the corresponding right Specht module, whose dimension is denoted by $f^\lambda=\dim_{\mathbb C}S^\lambda$. Let $\chi^\lambda$ be the character of $S^\lambda$, which is constant on each conjugacy class. Thus, we use
$\chi^\lambda(\mu)$ for the common value on $C_\mu$. 
For $\mu=(\mu_1,\ldots,\mu_t)\vdash r$, let
$a_j=1+\sum_{k<j}\mu_k$ for $1\le j\le t$, and define
\begin{equation}\label{eq:w-mu}
w_\mu=\prod_{j=1}^t
\bigl(s_{a_j}s_{a_j+1}\cdots s_{a_j+\mu_j-2}\bigr),
\end{equation}
where the $j$th factor is $1$ if $\mu_j=1$.

Let
$e_\lambda\in Z(\mathbb C\mathfrak S_r)$ be the primitive central
idempotent associated with $S^\lambda$. It is known that
\begin{align}\label{eq:class-idem}
e_\lambda
=\frac{f^\lambda}{r!}
  \sum_{\mu\vdash r}\chi^\lambda(\mu)c_\mu \quad\text{and}\quad
c_\mu
=\sum_{\lambda\vdash r}
  \frac{r!\,\chi^\lambda(\mu)}{z_\mu f^\lambda}e_\lambda.
\end{align}
Thus, the character table 
\begin{equation}\label{eq:chartable}
\mathsf{X}(r):=(\chi^\lambda(\mu))_{
\mu,\lambda\in\operatorname{Par}(r)}
\end{equation}
is, up to invertible diagonal
row and column factors, the transition matrix between the class-sum
basis and the primitive-central-idempotent basis of $Z(\mathbb C\mathfrak S_r)$.

\subsection{Symmetric functions}
Let $\mathrm{Sym}$ be the graded ring of symmetric functions over $\mathbb C$, and let $\mathrm{Sym}^r$ denote its homogeneous component of degree $r$. For $k\geq1$, let $p_k$ and $h_k$ denote
the $k$-th power-sum and $k$-th complete symmetric functions, respectively. For
$\mu=(\mu_1,\ldots,\mu_{\ell(\mu)})\vdash r$, write
$p_\mu=\prod_a p_{\mu_a}$. So $p_\mu p_\delta=p_{\mu\cup\delta}$ clearly. It is known that (cf.
\cite[Chapter I, (2.14')]{Mac95})
\begin{equation}\label{eq:h-p}
h_k=\sum_{\mu\vdash k}z_{\mu}^{-1}p_\mu.
\end{equation}

Let $s_\mu$ denote the Schur function corresponding to $\mu$. 
The well-known Frobenius formula \cite[Chapter~I, (7.8)]{Mac95} says that the character table 
$\mathsf{X}(r)$ is also the transition matrix between power-sum and Schur functions as follows:
\begin{equation}\label{eq:frob}
p_\mu=\sum_{\lambda\vdash r}\chi^\lambda(\mu)s_\lambda,\quad(\forall \mu\vdash r),
\end{equation}
which is equivalent to 
$$s_\lambda=\sum_{\mu\vdash r}\frac{1}{z_\mu}\chi^\lambda(\mu)p_\mu, \quad(\forall \lambda\vdash r).$$

For some special partitions, we have the following triangularity property with respect to the order $\prec$ defined in \eqref{def:ps-order} for the power-sum expansions of Schur functions, which is a direct consequence of the Jacobi–Trudi identity. We include a proof for completeness.

\begin{lem}\label{lem:JT-leading}
Let
$\gamma=(\gamma_1,\ldots,\gamma_{m+1})$ be a partition with
$\gamma_{m+1}\geq m+1$. 
The power-sum expansion of the Schur function $s_\gamma$ is
\begin{equation}\label{eq:specials-p}
s_\gamma
=
\frac{1}{\gamma_1\cdots\gamma_{m+1}}p_\gamma
+
\sum_{\eta\prec\gamma}
c_{\gamma\eta}p_\eta,
\quad
c_{\gamma\eta}\in\mathbb C.
\end{equation}
More precisely, every $\eta$ occurring in the sum satisfies either
$\ell(\eta)>m+1$, or
$\ell(\eta)=m+1$ and $\eta\triangleright\gamma$.
\end{lem}

\begin{proof}
Expand the Jacobi--Trudi determinant \cite[Chapter~I, (3.4)]{Mac95} as
$$
s_\gamma
=
\sum_{\sigma\in\mathfrak S_{m+1}}
\operatorname{sgn}(\sigma)
\prod_{i=1}^{m+1} h_{a_i},
\quad \text{where}\quad
a_i=\gamma_i-i+\sigma(i).
$$
For every $\sigma\in\mathfrak S_{m+1}$ and every $i$, one has
$a_i\geq\gamma_{m+1}-(m+1)+1\geq1$.
By \eqref{eq:h-p}, we see that, in the summand indexed by $\sigma$, the unique term
involving exactly $m+1$ power-sum factors is $\operatorname{sgn}(\sigma)
\frac{1}{a_1\cdots a_{m+1}}p_{\gamma_\sigma}$.

Let $\gamma_\sigma$ be the decreasing rearrangement of $(a_1,\ldots,a_{m+1})$.
By Lemma~\ref{lem:perm-major},
$\gamma_\sigma\unrhd\gamma$, with equality if and only if
$\sigma=\mathrm{id}$. Thus, when $\sigma\neq\mathrm{id}$, the
unique term involving exactly $m+1$ power-sum factors is indexed
by a partition strictly dominating $\gamma$.

For $\sigma=\mathrm{id}$, one has $a_i=\gamma_i$ for every $i$,
and the corresponding term is
$\prod_{i=1}^{m+1}\frac{p_{\gamma_i}}{\gamma_i}
=
\frac{1}{\gamma_1\cdots\gamma_{m+1}}p_\gamma$.
No other permutation contributes to $p_\gamma$, so its coefficient
is exactly $\frac{1}{\gamma_1\cdots\gamma_{m+1}}$. Every remaining
term is indexed by a partition $\eta$ satisfying either
$\ell(\eta)>m+1$, or
$\ell(\eta)=m+1$ and $\eta\triangleright\gamma$. In either case, we have $\eta\prec\gamma$.
\end{proof}

\subsection{Supersymmetric functions}
Let $x=(x_1,\ldots,x_m)$ and $y=(y_1,\ldots,y_n)$. Since the power sums $p_1,p_2,\ldots$ freely generate $\mathrm{Sym}$ as a polynomial algebra over $\mathbb C$, the assignments
$$
p_k\mapsto p_k(x)-p_k(y)
=
\sum_{i=1}^{m}x_i^k-\sum_{j=1}^{n}y_j^k,
\qquad k\geq1,
$$
determine a unique $\mathbb C$-algebra homomorphism
$$
\rho_{m,n}\colon
\mathrm{Sym}\longrightarrow
\mathbb C[x_1,\ldots,x_m,y_1,\ldots,y_n].
$$
With the standard grading on the target polynomial ring,
$\rho_{m,n}$ is a graded algebra homomorphism. Denote
$$
\Lambda_{m|n}:=\operatorname{im}\rho_{m,n}
\qquad \text{and} \qquad
\Lambda_{m|n}^r:=\rho_{m,n}(\mathrm{Sym}^r).
$$
Thus, $\Lambda_{m|n}^r$ is the homogeneous component of degree $r$
of $\Lambda_{m|n}$.
For $f\in\mathrm{Sym}$, we write $f(x|y):=\rho_{m,n}(f)$. More precisely, if $f=F(p_1,\ldots,p_N)$ for some polynomial $F$, then
$$
f(x|y)
=
F\bigl(
p_1(x)-p_1(y),\ldots,
p_N(x)-p_N(y)
\bigr).
$$
Applying $\rho_{m,n}$ to the power-sum and Schur functions yields the supersymmetric functions $p_\mu(x|y)$ and $s_\mu(x|y)$ for any partition $\mu$. Furthermore, the Schur functions admit the following hook Schur theorem.

\begin{thm}[{\cite[Lemma~6.4 and Corollary~6.5]{BR87};
\cite[Chapter~I, \S3, Example~24(b),(c)]{Mac95}}]\label{thm:hook-schur}
For every $\lambda\vdash r$, one has
$s_\lambda(x|y)=0$ if and only if
$\lambda\notin H_{m,n}(r)$. Moreover,
$\{s_\lambda(x|y)~|~\lambda\in H_{m,n}(r)\}$ is a
$\mathbb C$-basis of $\Lambda_{m|n}^r$. In particular,
$\dim_{\mathbb C}\Lambda_{m|n}^r=|H_{m,n}(r)|$.
\end{thm}

Using Theorem~\ref{thm:hook-schur}, we apply $\rho_{m,n}$ to \eqref{eq:frob}
to obtain
\begin{equation}\label{eq:super-frob}
p_\mu(x|y)
=\sum_{\lambda\in H_{m,n}(r)}
  \chi^\lambda(\mu)s_\lambda(x|y).
\end{equation}
Moreover, when $m\leq n$, the equation \eqref{eq:specials-p} is valid for any partition $\gamma=(\gamma_1,\ldots,\gamma_{m+1})$ with $\gamma_{m+1}>n$, which is non-hook. Thus, we have 
\begin{equation}\label{eq:0=s=p}
0=s_\gamma(x|y)
=
\frac{1}{\gamma_1\cdots\gamma_{m+1}}p_\gamma(x|y)
+
\sum_{\eta\prec\gamma}
c_{\gamma\eta}p_\eta(x|y).
\end{equation}

\begin{thm}\label{thm:gl-p-basis}
For all $m,n\geq0$ with $m+n>0$, the family
$\{p_\mu(x|y)~|~\mu\in B_{m,n}(r)\}$ forms a $\mathbb C$-basis of $\Lambda_{m|n}^r$.
\end{thm}

\begin{proof}
By definition, it is clear that
$p_k(y|x)=-p_k(x|y)$ for every $k\geq1$.
Thus $\Lambda_{m|n}$ and $\Lambda_{n|m}$, realized using
$(x|y)$ and $(y|x)$ in the same polynomial ring, coincide.
Moreover, $p_\mu(y|x)=(-1)^{\ell(\mu)}p_\mu(x|y)$.
Since $B_{m,n}(r)=B_{n,m}(r)$, it suffices to assume that
$m\leq n$, in which case $B_{m,n}(r)=H_{m,n}(r)$.

Since $\Lambda_{m|n}^r$ is the image of
$\mathrm{Sym}^r$ under $\rho_{m,n}$, the family
$\{p_\lambda(x|y)~|~\lambda\vdash r\}$ spans
$\Lambda_{m|n}^r$. It is therefore enough to express every
non-hook power sum as a linear combination of hook-indexed power
sums.

Take a non-hook partition $\lambda=(\lambda_1,\lambda_2,\cdots)\vdash r$ with $\lambda_{m+1}>n$. Write
$\gamma=(\lambda_1,\ldots,\lambda_{m+1})$ and
$\delta=(\lambda_{m+2},\lambda_{m+3},\ldots)$ so that $\lambda=\gamma\cup\delta$,
where $\delta$ is allowed to be empty and
$p_\varnothing(x|y)=1$. Since
$\lambda_{m+1}>n\geq m$, it follows from \eqref{eq:0=s=p} that
\begin{equation}\label{eq:p=sump}
p_\gamma(x|y)
=
\sum_{\nu\prec\gamma}
\kappa_{\gamma\nu}p_\nu(x|y),
\quad \text{with}\, 
\kappa_{\gamma\nu}\in\mathbb C.
\end{equation}
Note that for every $\nu\prec\gamma$ either
$\ell(\nu)>m+1$, in which case
$\ell(\nu\cup\delta)>\ell(\lambda)$, or
$\ell(\nu)=m+1$ and $\nu\triangleright\gamma$, in which case
Lemma~\ref{lem:union-dom} gives
$\nu\cup\delta\triangleright\lambda$. Therefore, we always have
$\nu\cup\delta\prec\lambda$. 
Multiply by $p_\delta(x|y)$ on the both sides of \eqref{eq:p=sump} and using
$p_\nu p_\delta=p_{\nu\cup\delta}$. 

After collecting equal terms, we obtain
$$
p_\lambda(x|y)
=
\sum_{\eta\prec\lambda}
\kappa_{\lambda\eta}p_\eta(x|y),
\quad \text{with}\,
\kappa_{\lambda\eta}\in\mathbb C.
$$

We prove by induction with respect to the order $\prec$ that $p_\lambda(x|y)\in \mathrm{Span}\{p_\mu(x|y)\mid\mu\in H_{m,n}(r)\}$ for every
$\lambda\vdash r$. This is immediate if $\lambda\in H_{m,n}(r)$.
If $\lambda\notin H_{m,n}(r)$, the preceding straightening relation
expresses $p_\lambda(x|y)$ as a linear combination of terms indexed
strictly below $\lambda$, all of which belong to $\mathrm{Span}\{p_\mu(x|y)\mid\mu\in H_{m,n}(r)\}$ by induction.
Hence $\{p_\mu(x|y)~|~\mu\in H_{m,n}(r)\}$ spans $\Lambda_{m|n}^r$.

Finally, Theorem~\ref{thm:hook-schur} gives
$\dim_{\mathbb C}\Lambda_{m|n}^r=|H_{m,n}(r)|$.
The family $\{p_\mu(x|y)~|~\mu\in H_{m,n}(r)\}$ has exactly this cardinality and is
therefore a basis.
\end{proof}

\subsection{Odd symmetric functions}
Let
$$\Gamma=\mathbb C[p_1,p_3,p_5,\ldots]\subseteq\mathrm{Sym}$$ be the ring of odd symmetric functions, 
and $\Gamma^r$ its homogeneous component of degree $r$.
For $\lambda\in\operatorname{SP}(r)$, let $Q_\lambda$ denote the
Schur $Q$-function (see
\cite[Chapter~III, \S8]{Mac95}). It is known (see
\cite[Chapter~III, (8.9)]{Mac95}) that the set
$
\{Q_\lambda~|~\lambda\in\operatorname{SP}(r)\}
$
forms a basis of $\Gamma^r$.

Denote $y=(y_1,\ldots,y_n)$ as before.
Let
$$
\rho_n:\mathrm{Sym}\rightarrow
\mathbb C[y_1,\ldots,y_n]^{\mathfrak S_n}
$$
be the specialization obtained by setting
$y_{n+1}=y_{n+2}=\cdots=0$, e.g. $p_\mu(y)=\rho_n(p_\mu)$.
We use the same notation $\rho_n$ for its
restriction to $\Gamma$, and denote
$$
\Gamma_n^r=\rho_n(\Gamma^r).
$$
One has $Q_\lambda(y)=0$ whenever $\ell(\lambda)>n$.

\begin{lem}[{\cite[Chapter~I, \S6]{Mac95}}]\label{lem:q-length-tri}
Let $\lambda\vdash r$ be a partition with $\ell(\lambda)\leq n$, and let
$m_d(\lambda)$ be the multiplicity of $d$ as a part of $\lambda$.
Then
$$
p_\lambda(y)
=
\left(\prod_{d\geq1}m_d(\lambda)!\right)m_\lambda(y)
+
\sum_{\substack{\nu\vdash r\\
                  \ell(\nu)<\ell(\lambda)}}
a_{\lambda\nu}m_\nu(y),
\quad
a_{\lambda\nu}\in\mathbb Z_{\geq0}.
$$
In particular, if $\lambda$ is strict, then the coefficient of
$m_\lambda(y)$ is $1$.
\end{lem}

\begin{thm}\label{thm:q-p-basis}
For all $n,r\geq 1$, the sets
\begin{equation*}
\{p_\mu(y)~|~\mu\in\operatorname{OP}_n^G(r)\}\quad \text{and}\quad
\{Q_\lambda(y)~|~\lambda\in\operatorname{SP}_n(r)\}
\end{equation*}
are two bases of $\Gamma_n^r$. Particularly, \begin{equation}\label{eq:Q-dim}
\dim\Gamma_n^r=|\operatorname{SP}_{n}(r)|=|\operatorname{OP}_n^G(r)|.
\end{equation} 
\end{thm}

\begin{proof}
Since $\{Q_\lambda~|~\lambda\in\operatorname{SP}(r)\}$
forms a basis of $\Gamma^r$, we know that the truncation $\{Q_\lambda(y)~|~\lambda\in\operatorname{SP}(r)\}$ spans $\Gamma_n^r$. So $\dim\Gamma_n^r\leq |\operatorname{SP}_{n}(r)|=|\operatorname{OP}_n^G(r)|$. So we only need to prove the linear independence of $\{p_\mu(y)~|~\mu\in\operatorname{OP}_n^G(r)\}$, which forces $\dim\Gamma_n^r=|\operatorname{SP}_{n}(r)|=|\operatorname{OP}_n^G(r)|$ and then indicates the two bases.

For $\lambda\in\operatorname{SP}_n(r)$, let
$\mu_\lambda=G^{-1}(\lambda)$.
For a symmetric polynomial $f$, write $[f:m_\tau(y)]$
for the coefficient of $m_\tau(y)$ in its monomial symmetric
expansion. Let
$$
A_n^Q(r)
=
\bigl([p_{\mu_\lambda}(y):m_\tau(y)]\bigr)_
{\lambda,\tau\in\operatorname{SP}_n(r)}.
$$

Work in the integral symmetric-polynomial ring and reduce
coefficients modulo $2$. In characteristic $2$, the Frobenius map
gives
$$
\bigl(p_j(y)\bigr)^{2^a}=p_{2^aj}(y).
$$
If
$\mu=(1^{m_1}3^{m_3}\cdots)\in\operatorname{OP}(r)$ and
$m_j=\sum_{a\geq0}\epsilon_{j,a}2^a$, then
\begin{equation}\label{eq:glaisher-cong}
p_\mu(y)
\equiv
\prod_{\substack{j\geq 1\text{ odd},\ a\geq 0\\
                         \epsilon_{j,a}=1}}
p_{2^aj}(y)
=
p_{G(\mu)}(y)
\pmod 2.
\end{equation}

Order $\operatorname{SP}_n(r)$ by decreasing length, refining this
to a total order arbitrarily. By
Lemma~\ref{lem:q-length-tri}, the matrix
$$
B_n^Q(r)
=
\bigl([p_\lambda(y):m_\tau(y)]\bigr)_
{\lambda,\tau\in\operatorname{SP}_n(r)}
$$
is triangular with every diagonal entry equal to $1$. Hence
$\det B_n^Q(r)=1$. Taking $\mu=\mu_\lambda$ in
\eqref{eq:glaisher-cong} gives
$A_n^Q(r)\equiv B_n^Q(r)\pmod2$. Thus $\det A_n^Q(r)$ is odd, and
in particular nonzero. Therefore, $\{p_\mu(y)~|~\mu\in\operatorname{OP}_n^G(r)\}$ is linearly independent, which implies $\dim\Gamma_n^r=|\operatorname{SP}_{n}(r)|=|\operatorname{OP}_n^G(r)|$. Thus $\{p_\mu(y)~|~\mu\in\operatorname{OP}_n^G(r)\}$ and 
$\{Q_\lambda(y)~|~\lambda\in\operatorname{SP}_n(r)\}$ are both bases of $\Gamma_n^r$.
\end{proof}

\section{Hecke and Hecke-Clifford algebras} \label{sec:quantum}

\subsection{Hecke algebras}

The Hecke algebra $\mathcal H_r$ over $\mathbb{K}$ is
generated by $T_1,\ldots,T_{r-1}$, subject to
\begin{equation}\label{rel:hecke}
(T_i-q)(T_i+1)=0, \quad T_iT_{i+1}T_i=T_{i+1}T_iT_{i+1},\quad
T_iT_j=T_jT_i \quad(|i-j|>1).
\end{equation}

The Hecke algebra $\mathcal H_r$ is split semisimple. Its
simple modules are the Specht modules $S_q^\lambda$ indexed by
$\lambda\vdash r$ (cf. \cite{DJ86}). We denote by $e_{\lambda,q}$
the corresponding primitive central idempotent and by
$\chi_q^\lambda$ the character of $S_q^\lambda$.

For $w\in\mathfrak S_r$ with a reduced form $w=s_{i_1}\cdots s_{i_{\ell(w)}}$, write
$T_w=T_{i_1}\cdots T_{i_{\ell(w)}}$, which is independent of the choice of a reduced form. The set $\{T_w~|~w\in\mathfrak S_r\}$ forms a basis of $\mathcal H_r$.

\subsection{Ordinary \texorpdfstring{$q$}{q}-Schur superalgebras}

Let $V=\mathbb C^{m|n}$ and $V_q=\mathbb{K}^{m|n}$ have homogeneous
bases, both denoted by $v_1,\ldots,v_{m+n}$, where the parity $|v_a|=\bar0$
for $a\leq m$ and $|v_a|=\bar1$ for $a>m$. For
$\mathbf i=(i_1,\ldots,i_r)\in\{1,\ldots,m+n\}^r$, let
$v_{\mathbf i}=v_{i_1}\otimes\cdots\otimes v_{i_r}$.
The right action of
$\mathcal H_r$ on $V_q^{\otimes r}$ is given, for $1\leq k<r$, by
\begin{equation}\label{eq:signed-q}
v_{\mathbf i}T_k=
\begin{cases}
(-1)^{|v_{i_k}||v_{i_{k+1}}|}
v_{\mathbf i s_k},
& i_k<i_{k+1},\\
qv_{\mathbf i},
& i_k=i_{k+1}\leq m,\\
-v_{\mathbf i},
& i_k=i_{k+1}>m,\\
q(-1)^{|v_{i_k}||v_{i_{k+1}}|}
v_{\mathbf i s_k}+(q-1)v_{\mathbf i},
& i_k>i_{k+1},
\end{cases}
\end{equation}
where $\mathbf i s_k$ exchanges $i_k$ and $i_{k+1}$.
These operators correspond to
$v\check R_k$ from \cite[(8.1.1), Lemma~8.2]{DR11} in the basis
$v^{-\operatorname{inv}(\mathbf i)}v_{\mathbf i}$, where
$\operatorname{inv}(\mathbf i)$ is the number of inversions of
$\mathbf i$. At $v=1$, they give the signed place-permutation
action of $\mathfrak S_r$ on $V^{\otimes r}$.



\begin{lem}\label{cor:q-hook-support}
For $\lambda\vdash r$, $V_q^{\otimes r}e_{\lambda,q}\neq0$ if and only if $\lambda\in H_{m,n}(r)$.
\end{lem}

\begin{proof}
See \cite[Corollary~7.8 and Proposition~8.3]{DR11}.
\end{proof}

%

The ordinary $q$-Schur superalgebra is defined as
$$
\mathcal S_q(m|n,r)=\operatorname{End}_{\mathcal H_r}(V_q^{\otimes r}).
$$

%

\subsection{Sergeev superalgebras and Hecke-Clifford algebras}
Let $\mathcal C_r$ be the Clifford superalgebra over $\mathbb{C}$ generated by odd
elements $c_1,\ldots,c_r$, subject to
\begin{equation}\label{rel:clifford}
c_i^2=1 \quad\text{and}\quad c_ic_j=-c_jc_i\quad (i\neq j).
\end{equation}
The Sergeev superalgebra is defined as
$\mathfrak S_r^c=\mathcal C_r\rtimes\mathbb C\mathfrak S_r$. Here, we collect some basic notation and facts for $\mathfrak S_r^c$ (cf. \cite[\S2.2]{WW13}). Simple $\mathfrak S_r^c$-supermodules are indexed by $\operatorname{SP}(r)$, up to parity shift.
For $\lambda\in\operatorname{SP}(r)$, let $U_1^\lambda$ be a
representative, which is of type
$\mathsf{M}$ if $\ell(\lambda)$ is even, and of type $\mathsf{Q}$ if
$\ell(\lambda)$ is odd. Let $e_\lambda^c$ be the primitive even
central idempotents corresponding to $\lambda\in\operatorname{SP}(r)$. 
Recall $\operatorname{SP}^{-}(r):=\{\lambda\in\operatorname{SP}(r)\mid
\ell(\lambda)\text{ is odd}\}$.
For $\lambda\in\operatorname{SP}^{-}(r)$, take
$\omega_\lambda\in Z(\mathfrak S_r^c e_\lambda^c)_{\bar1}$
the odd central element satisfying $\omega_\lambda^2=e_\lambda^c$.
Let $\zeta^\lambda$ be the ordinary character of $U_1^\lambda$. Since the restriction of
$\zeta^\lambda$ to $\mathbb C\mathfrak S_r$ is a class function, we denote by $\zeta^\lambda(\mu)$ the common value $\zeta^\lambda(w)$ for $w\in\mathfrak{S}_r$ with cycle type $\mu$.
The Frobenius formula for $\mathfrak S_r^c$ (see \cite[Theorem~3.7 and Remark~3.8]{WW13}) is
\begin{equation}\label{eq:spin-frob}
2^{\ell(\mu)}p_\mu
=\sum_{\lambda\in\operatorname{SP}(r)}
  2^{-\lceil\frac{\ell(\lambda)}{2}\rceil}
  \zeta^\lambda(\mu)Q_\lambda,
\quad
\mu\in\operatorname{OP}(r),
\end{equation}
where $\lceil x\rceil$ stands for the ceiling function (i.e., the least integer greater than or equal to $x\in\mathbb{R}$).
We denote by 
\begin{equation}\label{def:Yr}
\mathsf{Y}(r)=(\zeta^\lambda(\mu))_{ \mu\in\operatorname{OP}(r), \lambda\in\operatorname{SP}(r)}
\end{equation}
the character table.

The Hecke-Clifford algebra $\mathcal H_r^c$ over $\mathbb{K}$
is generated by even generators $T_1,\ldots,T_{r-1}$ and odd generators
$c_1,\ldots,c_r$, subject to
the Hecke relations \eqref{rel:hecke}, the Clifford relations \eqref{rel:clifford}, and their compatible conditions
$$
T_ic_j=c_jT_i \quad(j\neq i,i+1),\quad
T_ic_i=c_{i+1}T_i, \quad
T_ic_{i+1}=c_iT_i-(q-1)(c_i-c_{i+1}).
$$
This is the Hecke-Clifford algebra of \cite[\S2.3]{WW13}, with scalars
extended to $\mathbb K$, where the parameter $v$ of \cite{WW13} is our $q$.

Its simple supermodules $U_q^\lambda$, are also indexed by
$\lambda\in\operatorname{SP}(r)$; see \cite[Theorem~6.7 \& Corollary~6.8]{JN99}. Let $\zeta_q^\lambda$ be the
corresponding character, and $\zeta_q^\lambda(\mu)$ ($\mu\in\operatorname{OP}(r)$) its value following the character-table convention in \cite[\S4.2]{WW13}. It is known in \cite[Proposition~2.5]{WW13} that $\zeta_q^\lambda(\mu)|_{v=1}=\zeta^\lambda(\mu)$.

Since $\mathbb{K}=\overline{\mathbb{C}(v)}$ is algebraically closed, there also exist central idempotents and odd central elements in the split semisimple algebra $\mathcal H_{r}^c$, which we denote, respectively, by
$$
e_{\lambda,q}^c\in Z(\mathcal H_r^c)_{\bar{0}}
\quad(\lambda\in\operatorname{SP}(r))\quad\text{and}\quad \omega_{\lambda,q}\in Z(\mathcal H_r^c)_{\bar{1}} \quad(\lambda\in\operatorname{SP}^-(r)),
$$
with $\omega_{\lambda,q}^2=e_{\lambda,q}^c$.


\subsection{Queer \texorpdfstring{$q$}{q}-Schur superalgebras}
Assume $n\geq1$. Denote $W=\mathbb C^{n|n}$ and $W_q=\mathbb K^{n|n}$. Each has a homogeneous basis
$w_1,\ldots,w_n,\allowbreak w_{-1},\ldots,w_{-n}$,
where $w_a$ has parity $\widehat a=\bar{0}$ for $a>0$ and
$\widehat a=\bar{1}$ for $a<0$.
For $\mathbf i=(i_1,\ldots,i_r)\in\{\pm1,\ldots,\pm n\}^r$,
write $w_{\mathbf i}=w_{i_1}\otimes\cdots\otimes w_{i_r}$. The set $\{w_{\mathbf i}\mid\mathbf i\in\{\pm1,\ldots,\pm n\}^r\}$ forms a basis of $W_q^{\otimes r}$.
Let $\mathbf i s_k$ interchange $i_k$ and $i_{k+1}$, and let
$\mathbf i_k^-$ replace these entries by $-i_k$ and $-i_{k+1}$.

The Clifford generators act on the right by
$$
w_{\mathbf i}c_j
=\sqrt{-1}\,(-1)^{\sum_{a=1}^{j}\widehat{i_a}}\,
w_{(i_1,\ldots,i_{j-1},-i_j,i_{j+1},\ldots,i_r)},
\qquad 1\leq j\leq r.
$$
For $1\leq k<r$, the Hecke generators act on $W_q^{\otimes r}$ by
\begin{equation}\label{action:Hecke-Clifford}
w_{\mathbf i}T_k=
\begin{cases}
-w_{\mathbf i}-(q-1)w_{\mathbf i_k^-},
& i_k=i_{k+1}>0,\\[2pt]
q\,w_{\mathbf i},
& i_k=i_{k+1}<0,\\[2pt]
(q-1)w_{\mathbf i}-w_{\mathbf i s_k},
& i_k=-i_{k+1}>0,\\[2pt]
-q\,w_{\mathbf i s_k},
& i_k=-i_{k+1}<0,\\[2pt]
-v\,w_{\mathbf i s_k}-(q-1)w_{\mathbf i_k^-},
& |i_k|<|i_{k+1}|,\ i_{k+1}>0,\\[2pt]
(q-1)w_{\mathbf i}-(-1)^{\widehat{i_k}}v\,w_{\mathbf i s_k},
& |i_k|<|i_{k+1}|,\ i_{k+1}<0,\\[2pt]
-v\,w_{\mathbf i s_k}
+(q-1)\bigl(w_{\mathbf i}-(-1)^{\widehat{i_{k+1}}}w_{\mathbf i_k^-}\bigr),
& |i_k|>|i_{k+1}|,\ i_k>0,\\[2pt]
-(-1)^{\widehat{i_{k+1}}}v\,w_{\mathbf i s_k},
& |i_k|>|i_{k+1}|,\ i_k<0.
\end{cases}
\end{equation}
These formulas come from \cite[(6.1)--(6.2)]{DW18} via the
anti-involution $T_k\mapsto -qT_k^{-1}=q-1-T_k$, $c_j\mapsto c_j$ in
\cite[Lemma~2.1(3)]{DW18}, with each Clifford operator multiplied
by $\sqrt{-1}$. The action matrices have entries in
$\mathbb C[v,v^{-1}]$. 

We note that the actions $\eqref{action:Hecke-Clifford}$ are not compatible with those in \eqref{eq:signed-q} for $m=n$. From the Schur–Weyl duality viewpoint, this is reasonable because the quantum supergroup $\mathbf{U}_q(\mathfrak{q}_n)$ is not a subalgebra of the quantum supergroups $\mathbf{U}_q(\mathfrak{gl}_{n|n})$.

\begin{lem}\label{cor:q-strict-support}
For every $\lambda\in\operatorname{SP}(r)$, 
$W_q^{\otimes r}e_{\lambda,q}^c\neq0$ if and only if $\lambda\in\operatorname{SP}_n(r)$.
\end{lem}
\begin{proof}
See \cite[Corollary~6.3, Lemma~7.5 and Theorem~7.6]{DW18}.
\end{proof}

The queer $q$-Schur superalgebra is defined as
$$
\mathcal Q_q(n,r)
=\operatorname{End}_{\mathcal H_r^c}(W_q^{\otimes r}).
$$


\subsection{Submatrices of character tables}
\label{subsec:res-matrices}

Recall in \eqref{eq:chartable} the character table $\mathsf{X}(r)=(\chi^\lambda(\mu))_{\mu,\lambda\in\operatorname{Par}(r)}$ of $\mathfrak S_r$, with rows (resp. columns)
indexed by $\mu$ (resp. $\lambda$). Its submatrix
\begin{equation}\label{eq:gl-matrix}
\mathsf{X}_{m,n}(r)
=\bigl(\chi^\lambda(\mu)\bigr)_{\mu\in B_{m,n}(r),
             \lambda\in H_{m,n}(r)}
\end{equation}
is square, thanks to \eqref{eq:BH-card}. By
\eqref{eq:super-frob}, it is the coefficient matrix
from the selected supersymmetric power sums to the hook Schur basis.

\begin{lem}
\label{thm:gl-mat-ns}
For $r\geq 1$ and $m,n\geq0$ with $m+n>0$, the matrix $\mathsf X_{m,n}(r)$ is invertible.
\end{lem}

\begin{proof}
For each $\mu\in B_{m,n}(r)$, the restricted Frobenius formula
\eqref{eq:super-frob} expresses $p_\mu(x|y)$ in
the hook Schur basis $\{s_\lambda(x|y)\mid\lambda\in H_{m,n}(r)\}$, with coefficients
$\chi^\lambda(\mu)$, ($\lambda\in H_{m,n}(r)$).
By Theorem~\ref{thm:gl-p-basis}, the families
$\{p_\mu(x|y)\mid\mu\in B_{m,n}(r)\}$ is also a basis of
$\Lambda_{m|n}^{r}$. As the transition matrix for these two bases, $\mathsf X_{m,n}(r)$ must be invertible.
\end{proof}

Recall in \eqref{def:Yr} the character table $\mathsf Y(r)=(\zeta^\lambda(\mu))_{ \mu\in\operatorname{OP}(r), \lambda\in\operatorname{SP}(r)}$.
Its submatrix
\begin{equation}\label{eq:q-matrix-class}
\mathsf Y_n(r)
=
\bigl(\zeta^\lambda(\mu)\bigr)_{
\mu\in\operatorname{OP}_n^G(r),\ \lambda\in\operatorname{SP}_n(r)}
\end{equation}
is a square matrix by \eqref{eq:q-index-card}.

\begin{lem}\label{thm:q-mat-class}
For $n, r\geq 1$, the matrix $\mathsf Y_{n}(r)$ is invertible.
\end{lem}

\begin{proof}
Applying $\rho_n$ on \eqref{eq:spin-frob}, we have that, for
$\mu\in\operatorname{OP}_n^G(r)$,
$$
p_\mu(y)
=
\sum_{\lambda\in\operatorname{SP}_n(r)}
2^{-\ell(\mu)}\zeta^\lambda(\mu)
2^{-\lceil\frac{\ell(\lambda)}{2}\rceil}
Q_\lambda(y).
$$
Since $\{p_\mu(y)~|~\mu\in\operatorname{OP}_n^G(r)\}$ and $\{Q_\lambda(y)\mid\lambda\in\operatorname{SP}_n(r)\}$ are two bases of $\Gamma_n^r$ by Theorem~\ref{thm:q-p-basis}, the coefficient matrix
$$
\operatorname{diag}_{\mu\in\operatorname{OP}_n^G(r)}
  \bigl(2^{-\ell(\mu)}\bigr)
\mathsf Y_n(r)
\operatorname{diag}_{\lambda\in\operatorname{SP}_n(r)}
  \left(2^{-\lceil\frac{\ell(\lambda)}{2}\rceil}\right)
$$
is invertible. The two diagonal factors are invertible, so
$\mathsf Y_{n}(r)$ is invertible.
\end{proof}

\begin{rmk}
In fact, for every subset $S\subseteq\operatorname{SP}(r)$, the matrix
$$
\bigl(\zeta_q^\rho(G^{-1}(\lambda))\bigr)_{\lambda,\rho\in S}
$$
is invertible, both over $\mathbb{K}$ and at $v=1$.
In other words, once its rows are relabeled by the Glaisher
correspondence, every principal submatrix of the character table
of $\mathcal H_r^c$ is invertible.
The analogous statement for $\mathcal H_r$ follows over $\mathbb C(q)$
from the quantum Frobenius formula
\cite[Proposition~5.2]{WW15}, but fails in general at $v=1$;
for instance, $\chi^{(3,1)}((3,1))=0$.
\end{rmk}

\section{Centers} \label{sec:center}
In this section, we provide an approach for studying centers of endomorphism superalgebras, which also applies to ordinary algebras regarded as superalgebras with trivial odd parts. 

\subsection{Centers of endomorphism superalgebras}
Let $F$ be a field with $\operatorname{char}F\neq2$. Let $A=A_{\bar{0}}\oplus A_{\bar{1}}$ be a finite-dimensional split semisimple $F$-superalgebra. That is, $A$ is isomorphic to a direct sum of split simple $F$-superalgebras, each of type $\mathsf{M}$ or type $\mathsf{Q}$.
Let $\{e_\lambda\mid\lambda\in\Lambda\}$ be the set of
primitive even central idempotents of $A$.
For each $\lambda\in\Lambda$, choose a simple right
$A$-supermodule $L^\lambda$ such that
$L^\lambda e_\lambda=L^\lambda$.
Let
$$
\Lambda_{\bar0}=\Lambda,
\qquad
\Lambda_{\bar1}
=\{\lambda\in\Lambda\mid Ae_\lambda\text{ is of type $\mathsf{Q}$}\}.
$$
For $\lambda\in \Lambda$ such that $Ae_\lambda$ is of type $\mathsf{Q}$, choose $\omega_\lambda\in Z(Ae_\lambda)_{\bar 1}$ with $\omega_\lambda^2=e_\lambda$.
Write
$$
u_{\lambda,\bar0}=e_\lambda
\quad(\lambda\in\Lambda_{\bar0}),
\qquad
u_{\lambda,\bar1}=\omega_\lambda
\quad(\lambda\in\Lambda_{\bar1}).
$$
For each $\bar i\in\mathbb Z_2$, the family
$\{u_{\lambda,\bar i}\mid\lambda\in\Lambda_{\bar i}\}$
is an $F$-basis of $Z(A)_{\bar i}$.

For a nonzero finite-dimensional right $A$-supermodule $M$, let
$\eta_M:A^{\operatorname{op}}\to\End_F(M)$
be the superalgebra homomorphism given by $\eta_M(a)(v)=va$, where $A^{\operatorname{op}}$ is the
ordinary opposite algebra. Let $B=\End_A(M)$ be the
superalgebra of all $A$-linear endomorphisms of $M$, with its natural
grading and multiplication given by composition. Both even and
odd maps satisfy $f(va)=f(v)a$. Denote by
$\Lambda_M=\{\lambda\in\Lambda~|~Me_\lambda\neq0\}$ the support of $M$, and let $\Lambda_{M,\bar i}=\Lambda_M\cap\Lambda_{\bar i}$ for $\bar i\in \mathbb Z_2$.

\begin{thm}\label{prop:supp-center}
The homomorphism $\eta_M$ induces a surjective graded algebra homomorphism $Z(A)\to Z(B)$, with 
$$
\ker(\eta_M|_{Z(A)})=\bigoplus\limits_{\lambda\notin\Lambda_M}Z(Ae_\lambda).
$$
For each $\bar i\in\mathbb Z_2$, we have
$$
Z(B)_{\bar i}=\bigoplus_{\lambda\in\Lambda_{M,\bar i}}F\eta_M(u_{\lambda,\bar i})
$$
and hence $\dim_F Z(B)_{\bar i}=|\Lambda_{M,\bar i}|$.
The elements $\eta_M(e_\lambda)$, $\lambda\in\Lambda_M$,
are the primitive even central idempotents of $B$.
\end{thm}

\begin{proof}
For homogeneous $z\in Z(A)$, the right multiplication by $z$ has parity
$|z|$ and is $A$-linear, since $(va)z=(vz)a$ for $v\in M$ and
$a\in A$. For every $f\in B$ and $v\in M$, one has
$f(vz)=f(v)z$, so $\eta_M(z)\in Z(B)_{|z|}$.
Moreover,
$\eta_M(z)\eta_M(z')=\eta_M(z'z)=\eta_M(zz')$
for $z,z'\in Z(A)$. Hence $\eta_M$ restricts to a graded algebra
homomorphism $Z(A)\to Z(B)$.

The decomposition
$M=\bigoplus_{\lambda\in\Lambda_M}Me_\lambda$ gives
$$
B\cong
\bigoplus_{\lambda\in\Lambda_M}
B_\lambda,\qquad \text{where} \quad B_\lambda=\End_A(Me_\lambda).
$$
Each $B_\lambda$ is split simple with identity $\eta_M(e_\lambda)$. If $Ae_\lambda$ is of type $\mathsf{M}$, both centers are scalar and are identified by $\eta_M$. If it is of type $\mathsf{Q}$,  then $\eta_M(\omega_\lambda)^2=\eta_M(e_\lambda)\ne0$. Since $Z(B_\lambda)_{\bar 1} $ has dimension at most one, it is spanned by $\eta_{M}(\omega_\lambda)$. Hence $\eta_M$ induces an isomorphism $Z(Ae_\lambda)\to Z(B_\lambda)$ for every $\lambda\in \Lambda_M$.

The blockwise isomorphisms, together with $Me_\lambda=0$ for $\lambda\notin\Lambda_M$, give surjectivity, the kernel formula and the two center decompositions. These decompositions yield the dimension formulas and show that the $\eta_M(e_\lambda)$, $\lambda\in\Lambda_M$, are precisely the primitive even central idempotents of $B$.
\end{proof}

Fix $\bar i\in\mathbb Z_2$.
Let $(\xi_\mu)_{\mu\in\Omega}$ be a finite family of
elements of $Z(A)_{\bar i}$, and write
\begin{equation}\label{eq:center-expand}
\xi_\mu
=
\sum_{\lambda\in\Lambda_{\bar i}}
c_{\mu\lambda}u_{\lambda,\bar i},
\qquad \text{with}\quad c_{\mu\lambda}\in F.
\end{equation}
Denote $\widehat{\xi}_\mu=\eta_M(\xi_\mu)$.
For $I\subseteq\Omega$, let
$C_{I,\Lambda_{M,\bar i}}
=(c_{\mu\lambda})_{\mu\in I,\ \lambda\in\Lambda_{M,\bar i}}$
be the restricted coefficient matrix.

\begin{prop}\label{prop:coeff-matrix}
For $\mu\in\Omega$, we have
$$
\widehat{\xi}_\mu
=\sum_{\lambda\in\Lambda_{M,\bar i}}
c_{\mu\lambda}\eta_M(u_{\lambda,\bar{i}}).
$$
Moreover, for every $I\subseteq\Omega$,
$$
\dim_F\Span_F\{\widehat{\xi}_\mu\mid\mu\in I\}
=\rank_F C_{I,\Lambda_{M,\bar i}}.
$$
In particular, the family
$(\widehat{\xi}_\mu)_{\mu\in I}$ is an $F$-basis of
$Z(B)_{\bar i}$ if and only if
$C_{I,\Lambda_{M,\bar i}}$ is an invertible square matrix.
\end{prop}

\begin{proof}
By Theorem~\ref{prop:supp-center},
$\eta_M(u_{\lambda,\bar i})=0$ for
$\lambda\in\Lambda_{\bar i}\setminus\Lambda_{M,\bar i}$,
and
$\{\eta_M(u_{\lambda,\bar i})\mid\lambda\in\Lambda_{M,\bar i}\}$
is an $F$-basis of $Z(B)_{\bar i}$.
Applying $\eta_M$ to \eqref{eq:center-expand} gives the
displayed expansion.
The rows of $C_{I,\Lambda_{M,\bar i}}$ are therefore the
coordinates of the elements $\widehat{\xi}_\mu$, $\mu\in I$,
in this basis.
The rank formula and the basis criterion follow.
\end{proof}

\subsection{Casimir elements of symmetric superalgebras}
Assume further that $A$ is symmetric, equipped with an even symmetric trace form $\tau:A\to F$. It means that 
$\tau(ab)=\tau(ba)$ for all $a,b\in A$, the trace $\tau$ vanishes
on $A_{\bar1}$, and the bilinear form
$(a,b)\mapsto\tau(ab)$ is nondegenerate. 

Let $\mathcal B$ be a homogeneous basis of $A$, and
$\mathcal B^\vee=\{b^\vee\mid b\in\mathcal B\}$ its
$\tau$-dual basis, characterized by
$\tau(b^\vee b')=\delta_{b,b'}$.
The Casimir element associated with $h\in A$ is defined as
\begin{equation}\label{eq:casimir-map}
C_\tau(h)=\sum_{b\in\mathcal B}b^\vee h b.
\end{equation}

\begin{prop}\label{prop:casimir-form}
The Casimir element $C_\tau(h)$ is independent of the choice of $\mathcal B$,
and $C_\tau(A_{\bar i})\subseteq Z(A)_{\bar i}$ for
$\bar i\in\mathbb Z_2$. There exist scalars
$\kappa_\lambda\in F^\times$ such that, $\forall h\in A_{\bar{i}}$, $\bar i\in\mathbb{Z}_2$,
\begin{equation}\label{eq:casimir-expand}
C_\tau(h)
=\sum_{\lambda\in\Lambda_{\bar i}}
  \kappa_\lambda\chi^\lambda(h u_{\lambda,\bar i})u_{\lambda,\bar i},
\end{equation}
where $\chi^\lambda$ is the ordinary character of $L^\lambda$.
\end{prop}

\begin{proof}
It is a standard result for dual bases that the tensor
$\sum_{b\in\mathcal B}b^\vee\otimes b$ is independent of
$\mathcal B$ and satisfies
$$
\sum_{b\in\mathcal B}ab^\vee\otimes b
=
\sum_{b\in\mathcal B}b^\vee\otimes ba
\qquad(a\in A).
$$
Applying $x\otimes y\mapsto xhy$ proves that $C_\tau(h)$ is central
and independent of $\mathcal B$. Since $|b^\vee|=|b|$ and $\tau$ is even, the map $C_\tau$ preserves parity.

For $h\in A_{\bar0}$, regard $L^\lambda$ as a left
$A^{\operatorname{op}}$-module via $\rho_\lambda(a)(v)=va$.
The same trace $\tau$ symmetrizes $A^{\operatorname{op}}$, and
$\mathcal B$ is dual to $\mathcal B^\vee$ there.
Using these bases, \cite[Lemma~5.2 \& Proposition~5.4]{WW13} gives
$$
\rho_\lambda(C_\tau(h))
=
\sum_{b\in\mathcal B}
\rho_\lambda(b)\rho_\lambda(h)\rho_\lambda(b^\vee)
=
\kappa_\lambda\chi^\lambda(h)\id_{L^\lambda},
$$
where $\kappa_\lambda\in F^\times$ is independent of $h$.
Since the direct sum of the $L^\lambda$ is faithful and
$\chi^\lambda(he_\lambda)=\chi^\lambda(h)$, this proves
\eqref{eq:casimir-expand} for $\bar i=\bar0$.

For $h\in A_{\bar1}$ and $\lambda\in\Lambda_{\bar1}$,
the element $h\omega_\lambda$ lies in $(Ae_\lambda)_{\bar0}$.
Since $\omega_\lambda\in Z(A)$, the definition of $C_\tau$
gives $C_\tau(h\omega_\lambda)=C_\tau(h)\omega_\lambda$.
Using the even case and $\omega_\lambda^2=e_\lambda$, we obtain
$$
C_\tau(h)e_\lambda
=
C_\tau(h\omega_\lambda)\omega_\lambda
=
\kappa_\lambda\chi^\lambda(h\omega_\lambda)\omega_\lambda.
$$
Summing over $\lambda\in\Lambda_{\bar1}$ proves
\eqref{eq:casimir-expand} for $\bar i=\bar1$.
\end{proof}

For a finite family $(h_\mu)_{\mu\in I}$ in $A_{\bar i}$ with
$|I|=|\Lambda_{M,\bar i}|$, Propositions~\ref{prop:coeff-matrix}
\& \ref{prop:casimir-form} show that
$(\eta_M(C_\tau(h_\mu)))_{\mu\in I}$ is a basis of $Z(B)_{\bar i}$
if and only if the matrix
$$
\bigl(\chi^\lambda(h_\mu u_{\lambda,\bar i})\bigr)_
{\mu\in I,\ \lambda\in\Lambda_{M,\bar i}}
$$
is invertible.

\subsection{Centers of Hecke algebras}
In fact, the Hecke algebra $\mathcal H_r$ is a symmetric algebra with a canonical symmetrizing trace 
$$\tau:\mathcal H_r\to \mathbb{K}, \quad T_w\mapsto \delta_{w,\mathbbm{1}},$$
which induces a bilinear form 
$(T_w,T_\sigma)=\delta_{w,\sigma^{-1}}q^{\ell(w)}$.
Denote 
$$T_w^\vee=q^{-\ell(w)}T_{w^{-1}}
\quad(w\in\mathfrak S_r).
$$
Then $\{T_w^\vee\mid w\in\mathfrak S_r\}$ forms a basis of $\mathcal H_r$ dual to
$\{T_w\mid w\in\mathfrak S_r\}$ (cf. \cite[\S2.2]{WW15}).

For $\mu\vdash r$, denote
\begin{equation}\label{eq:gl-casimir}
\mathsf A_{\mu,q}
:=\sum_{w\in\mathfrak S_r}
T_w^\vee T_{w_\mu}T_w
=\sum_{w\in\mathfrak S_r}
q^{-\ell(w)}T_{w^{-1}}T_{w_\mu}T_w.
\end{equation}
Proposition~\ref{prop:casimir-form} gives
$\mathsf A_{\mu,q}\in Z(\mathcal H_r)$.
Moreover,
$$
\mathsf A_{\mu,q}\big|_{v=1}
=\sum_{w\in\mathfrak S_r}w^{-1}w_\mu w
=z_\mu c_\mu.
$$

\begin{lem}\label{lem:hecke-center-basis}
The set
$\{\mathsf A_{\mu,q}\mid\mu\vdash r\}$
forms a basis of $Z(\mathcal H_r)$.
\end{lem}

\begin{proof}
The coordinate matrix of this family in the standard basis
$\{T_w\mid w\in\mathfrak S_r\}$ has entries in $\mathbb C[v,v^{-1}]$.
Its specialization has full column rank, since
$\mathsf A_{\mu,q}\big|_{v=1}=z_\mu c_\mu$ ($\mu\vdash r$) form a basis
of $Z(\mathbb C\mathfrak S_r)$.
Thus $\{\mathsf A_{\mu,q}\mid\mu\vdash r\}$ is
linearly independent.
It must be a basis because
$\dim_{\mathbb K} Z(\mathcal H_r)=|\operatorname{Par}(r)|$.
\end{proof}

\subsection{Center of the ordinary \texorpdfstring{$q$}{q}-Schur superalgebra}
By Lemma~\ref{cor:q-hook-support} and
Theorem~\ref{prop:supp-center}, the tensor-space action
induces a surjective homomorphism
$$
Z(\mathcal H_r)
\xrightarrow{\ \eta_{q,r}\ }
Z(\mathcal S_q(m|n,r)),
$$
with
$$
\dim_{\mathbb K} Z(\mathcal S_q(m|n,r))
=
|H_{m,n}(r)|.
$$

For $\lambda\in H_{m,n}(r)$, let
$\varepsilon_{\lambda,q}=\eta_{q,r}(e_{\lambda,q})$.
These are the primitive central idempotents of $\mathcal S_q(m|n,r)$.

\begin{thm}\label{thm:gl-center-bases}
The set
$\{\eta_{q,r}(\mathsf A_{\mu,q})\mid\mu\in B_{m,n}(r)\}$
forms a basis of $Z(\mathcal S_q(m|n,r))$.
\end{thm}

\begin{proof}
By Proposition~\ref{prop:coeff-matrix} and Proposition~\ref{prop:casimir-form},
$$
\eta_{q,r}(\mathsf A_{\mu,q})
=
\sum_{\lambda\in H_{m,n}(r)}
s_{\lambda,q}\chi_q^\lambda(T_{w_\mu})\varepsilon_{\lambda,q},
$$
where $s_{\lambda,q}\in \mathbb K^\times$ is the Schur element associated
with $S_q^\lambda$ and $\tau$ (cf. \cite[\S2.3]{WW15}).
Since these scalars are nonzero, Proposition~\ref{prop:coeff-matrix}
reduces the assertion to the invertibility over $\mathbb K$ of the matrix
$$
\bigl(\chi_q^\lambda(T_{w_\mu})\bigr)_{
\mu\in B_{m,n}(r),\ \lambda\in H_{m,n}(r)}.
$$
By \cite[\S2.1]{WW15}, this matrix specializes at $v=1$ to
$\mathsf X_{m,n}(r)$, which is invertible by
Lemma~\ref{thm:gl-mat-ns}.
Its determinant is therefore nonzero in $\mathbb{K}$.
\end{proof}

\begin{rmk}\label{rmk:gl-other-bases}
The above argument also works for the Geck--Rouquier basis
and the Jones basis of $Z(\mathcal{H}_r)$. Let $\Gamma_\mu$ and $J_\mu$ be the
corresponding central elements (cf. \cite{GR97,Jo90}), normalized as in
\cite[\S2.2 and \S4.1]{WW15}, so that
$$
\Gamma_\mu|_{v=1}=c_\mu,
\qquad
J_\mu|_{v=1}=d_\mu c_\mu,$$
where
$d_\mu=\prod_{j\geq1}m_j(\mu)!\neq 0$.
Consequently, the sets
$$
\{\eta_{q,r}(\Gamma_\mu)\mid\mu\in B_{m,n}(r)\}
\quad\text{and}\quad
\{\eta_{q,r}(J_\mu)\mid\mu\in B_{m,n}(r)\}
$$
are both bases of $Z(\mathcal S_q(m|n,r))$. Very recently, this is done by Fu, Luo and Sun in \cite{FLS26}. 
\end{rmk}

\subsection{Centers of Hecke-Clifford algebras}
Let 
$$
\mathcal B_q
=
\{T_wc^\alpha~|~w\in\mathfrak S_r,\ \alpha\in\{0,1\}^r\},
\quad\text{where}\quad c^\alpha=c_1^{\alpha_1}\cdots c_r^{\alpha_r},
$$
be the standard homogeneous basis of $\mathcal H_r^c$.
Its specialization at $v=1$ is a basis
$\mathcal B_1=\{wc^\alpha\mid w\in\mathfrak S_r,\
\alpha\in\{0,1\}^r\}$ of $\mathfrak S_r^c$.

Let $\tau_q$ be the even symmetrizing trace on
$\mathcal H_r^c$ constructed in
\cite[\S5.2, Theorem~5.8]{WW13}.
Its values on $\mathcal B_q$ lie in $\mathbb C[v,v^{-1}]$
(see \cite[Theorem~4.7]{WW13}), and
$$
\tau_q(T_{w_\nu})
=\left(\frac{q-1}{2}\right)^{r-\ell(\nu)}
\qquad(\nu\in\operatorname{OP}(r)).
$$
At $v=1$, it specializes to the coefficient-of-$1$ trace
$\tau_1$ on $\mathfrak S_r^c$ with respect to $\mathcal B_1$;
see \cite[Remark~5.11]{WW13}.
Let $\mathcal B_q^\vee=\{b^\vee\mid b\in\mathcal B_q\}$
be the $\tau_q$-dual basis. As defined in \eqref{eq:casimir-map}, we obtain Casimir elements $C_{\tau_q}(h)=\sum_{b\in\mathcal B_q}b^\vee h b\in Z(\mathcal H_r^c)$ for $h\in\mathcal H_r^c$, and $C_{\tau_1}(h)=\sum_{b\in\mathcal B_1}b^\vee h b\in Z(\mathfrak S_r^c)$ for $h\in\mathfrak S_r^c$.

\begin{lem}\label{lem:q-casimir-special}
The coordinates of $\mathcal B_q^\vee$ in $\mathcal B_q$ are regular at $v=1$.
This $\tau_q$-dual basis specializes to the $\tau_1$-dual basis
$\mathcal B_1^\vee$.
If $h\in(\mathcal H_r^c)_{\bar i}$ has coordinates regular at $v=1$ with respect to $\mathcal B_q$, then
$C_{\tau_q}(h)\in Z(\mathcal H_r^c)_{\bar i}$
has the same property, and
$$
C_{\tau_q}(h)\big|_{v=1}
=
C_{\tau_1}\bigl(h\big|_{v=1}\bigr).
$$
\end{lem}

\begin{proof}
Let $M(q)=(\tau_q(bb'))_{b,b'\in\mathcal B_q}$ and
$M(1)=(\tau_1(bb'))_{b,b'\in\mathcal B_1}$. By the definition of dual bases,
$M(q)^{-1}$ is the coordinate matrix of $\mathcal B_q^\vee$ in $\mathcal B_q$,
and similarly at $v=1$.
The entries of $M(q)$ lie in $\mathbb C[v,v^{-1}]$, and $M(1)$ is
invertible, since $\tau_1(b^{-1}b')=\delta_{b,b'}$
for $b,b'\in\mathcal B_1$.
Thus $M(q)^{-1}$ is regular at $v=1$ and specializes to
$M(1)^{-1}$, which proves the assertions about the dual basis.
The remaining assertions follow from \eqref{eq:casimir-map}
and Proposition~\ref{prop:casimir-form}.
\end{proof}

For $\mu\in\operatorname{OP}(r)$, write
$\mathsf{C}_{\mu,1}=C_{\tau_1}(w_\mu)$ and
$\mathsf{C}_{\mu,q}=C_{\tau_q}(T_{w_\mu})$. For $\lambda=(\lambda_1,\ldots,\lambda_\ell)
\in\operatorname{SP}^{-}(r)$, let
$a_j=1+\sum_{k<j}\lambda_k$ for $1\leq j\leq\ell$, and define
$$
h_{\lambda,1}=c_{a_1}\cdots c_{a_\ell}w_\lambda,
\qquad
h_{\lambda,q}=c_{a_1}\cdots c_{a_\ell}T_{w_\lambda}.
$$
These elements are odd since $\ell$ is odd.
Let
\begin{equation}\label{eq:q-odd-class-casimir}
\begin{aligned}
\mathsf D_{\lambda,1}
&=C_{\tau_1}(h_{\lambda,1})
=\sum_{b\in\mathcal B_1}b^{-1}h_{\lambda,1}b,\\
\mathsf D_{\lambda,q}
&=C_{\tau_q}(h_{\lambda,q})
=\sum_{b\in\mathcal B_q}b^\vee h_{\lambda,q}b.
\end{aligned}
\end{equation} By Proposition~\ref{prop:casimir-form} and
Lemma~\ref{lem:q-casimir-special},
$\mathsf C_{\mu,q}\in Z(\mathcal H_{r}^c)_{\bar0}$
and
$\mathsf D_{\lambda,q}\in Z(\mathcal H_{r}^c)_{\bar1}$,
with
$$
\mathsf C_{\mu,q}\big|_{v=1}=\mathsf C_{\mu,1},
\qquad
\mathsf D_{\lambda,q}\big|_{v=1}=\mathsf D_{\lambda,1}.
$$

For $\rho\in\operatorname{SP}^{-}(r)$, write $\zeta^\rho=\zeta^{\rho,+}+\zeta^{\rho,-}$, where $\zeta^{\rho,+}$ and $\zeta^{\rho,-}$ are the ordinary irreducible characters afforded by the $(+1)$- and
$(-1)$-eigenspaces of $\omega_\rho$ in $U_1^\rho$,
respectively.
For $h\in\mathfrak S_r^c$, we have
$$
\zeta^\rho(h\omega_\rho)
=
\zeta^{\rho,+}(h)-\zeta^{\rho,-}(h).
$$
If $h$ is odd, then $\zeta^\rho(h)=0$, and hence
$\zeta^\rho(h\omega_\rho)=2\zeta^{\rho,+}(h)$.

\begin{lem}\label{lem:q-odd-values}
For $\lambda,\rho\in\mathrm{SP}^-(r)$, $\zeta^\rho(h_{\lambda,1}\,\omega_\rho)\neq0$ if and only if $\rho=\lambda$.
\end{lem}

\begin{proof}
Suppose first that $r\leq3$. Then $\operatorname{SP}^{-}(r)=\{(r)\}$, so we only need
to show that $\zeta^{(r)}(h\,\omega_{(r)})\neq0$, where $h=h_{(r),1}$.
A direct computation shows that exactly $2r$ of the conjugates $b^{-1}hb$
$(b\in\mathcal B_1)$ are equal to $h$, and none is equal to $-h$.
Since each conjugate lies in $\pm\mathcal B_1$, it follows that
$C_{\tau_1}(h)\neq0$. On the other hand, since $(r)$ is the only element of
$\operatorname{SP}^{-}(r)$, Proposition~\ref{prop:casimir-form} gives
$C_{\tau_1}(h)=\kappa\,\zeta^{(r)}(h\,\omega_{(r)})\,\omega_{(r)}$ for some
$\kappa\in\mathbb C^\times$, and the claim follows.

Now let $r\geq4$, and let $\tilde H_r$ be the double cover of the hyperoctahedral
group defined in \cite[(2.2)--(2.5)]{HJ17}, with generators $s_i$, $\tilde c_j$ and
central element $z$. By the defining relations, $s_i\mapsto s_i$ and
$\tilde c_j\mapsto\sqrt{-1}\,c_j$ induce an isomorphism
$\mathbb C\tilde H_r/(1+z)\cong\mathfrak S_r^c$.
Under this isomorphism, $h_{\lambda,1}$ corresponds, up to the scalar
$(\sqrt{-1})^{\ell(\lambda)}$, to an element of one of the odd split classes
$D^\pm_\lambda$ of \cite[\S3]{HJ17}. Moreover, the supercharacter $\xi_\nu$ of
\cite{HJ17} corresponds to $\zeta^\nu$, as one sees by comparing
\cite[(4.2)]{HJ17} with \eqref{eq:spin-frob}; hence $\xi^\pm_\nu$ corresponds to
$\zeta^{\nu,\pm}$ up to order. Since
$\zeta^\rho(h_{\lambda,1}\,\omega_\rho)=2\,\zeta^{\rho,+}(h_{\lambda,1})$, the claim
follows from \cite[Theorem~4.5]{HJ17}.
\end{proof}

Thanks to Lemma~\ref{lem:q-odd-values}, we have 
\begin{equation}\label{eq:Dv=1}
\mathsf D_{\lambda,q}|_{v=1}=2\kappa_{\lambda,1}
\zeta^{\lambda,+}(h_{\lambda,1})\omega_\lambda
\neq0,
\end{equation}
where $\kappa_{\lambda,1}\in\mathbb C^\times$ is the scalar in
Proposition~\ref{prop:casimir-form} for $(\mathfrak S_r^c,\tau_1)$. Particularly,
$\mathsf D_{\lambda,q}\neq 0$.

\begin{lem}\label{thm:hc-center-bases}
The sets
$\{\mathsf C_{\mu,q}\mid\mu\in\operatorname{OP}(r)\}$ and $\{\mathsf D_{\lambda,q}\mid\lambda\in\operatorname{SP}^{-}(r)\}$
form bases of $Z(\mathcal H_r^c)_{\bar{0}}$ and $Z(\mathcal H_r^c)_{\bar{1}}$, respectively.
\end{lem}

\begin{proof}
By Proposition~\ref{prop:casimir-form}, the coefficient matrix
of the even family relative to the primitive even central
idempotents is
$$
\bigl(\zeta_q^\lambda(\mu)\bigr)_{
\mu\in G^{-1}(\operatorname{SP}(r)),\ \lambda\in\operatorname{SP}(r)}
\operatorname{diag}_{\lambda\in\operatorname{SP}(r)}
\bigl(\kappa_{\lambda,q}\bigr),
$$
where the matrix $\bigl(\zeta_q^\lambda(\mu)\bigr)_{
\mu\in G^{-1}(\operatorname{SP}(r)),\ \lambda\in\operatorname{SP}(r)}$ is the transpose of the
character table of $\mathcal H_r^c$ and is invertible
by \cite[\S4.2]{WW13}.
Since every $\kappa_{\lambda,q}$ is nonzero, the coefficient
matrix is invertible over $\mathbb K$, and hence the even family is a basis of $Z(\mathcal H_r^c)_{\bar0}$.

Since $\{\omega_\lambda \mid \lambda\in\operatorname{SP}^-(r)\}$ forms a basis of $Z(\mathfrak S_r^c)_{\bar{1}}$, it follows from \eqref{eq:Dv=1} that $\mathsf D_{\lambda,q} (\lambda\in\operatorname{SP}^-(r))$ is linearly independent. Therefore, these $\mathsf D_{\lambda,q}$ must be a basis of $Z(\mathcal H_r^c)_{\bar1}$ because $\dim_{\mathbb K}Z(\mathcal H_r^c)_{\bar1}
=|\operatorname{SP}^{-}(r)|$.
\end{proof}

\subsection{Centers of queer \texorpdfstring{$q$}{q}-Schur superalgebras}

By Theorem~\ref{prop:supp-center}, the tensor-space action
induces a surjective graded algebra homomorphism
$$
Z(\mathcal H_r^c)
\xrightarrow{\ \eta_q^Q\ }
Z(\mathcal Q_q(n,r)).
$$
Thanks to Lemma~\ref{cor:q-strict-support}, the occurring blocks are indexed by $\operatorname{SP}_n(r)$, while the occurring type $\mathsf{Q}$ blocks are indexed by $\operatorname{SP}_n^{-}(r)$. Hence
$$
\dim_{\mathbb K}Z(\mathcal Q_q(n,r))_{\bar0}
=
|\operatorname{SP}_n(r)| \quad \text{and}\quad
\dim_{\mathbb K}Z(\mathcal Q_q(n,r))_{\bar1}
=
|\operatorname{SP}_n^{-}(r)|.
$$
For $\lambda\in\operatorname{SP}_n(r)$, let
$\varepsilon_{\lambda,q}^Q=\eta_q^Q(e_{\lambda,q}^c)$.
Then $\{\varepsilon_{\lambda,q}^Q \mid \lambda\in\operatorname{SP}_n(r)\}$ is the set of primitive even central idempotents of $\mathcal Q_q(n,r)$.

\begin{thm}\label{thm:q-even-bases}
The sets
$\{\eta_q^Q(\mathsf C_{\mu,q})
\mid\mu\in\operatorname{OP}_n^G(r)\}$ and $\{\eta_q^Q(\mathsf D_{\lambda,q})
\mid\lambda\in\operatorname{SP}_n^{-}(r)\}$
form bases of $Z(\mathcal Q_q(n,r))_{\bar{0}}$ and $Z(\mathcal Q_q(n,r))_{\bar{1}}$, respectively. 
\end{thm}

\begin{proof}
By Proposition~\ref{prop:coeff-matrix} and Proposition~\ref{prop:casimir-form},
the coefficient matrix of the even family $\{\eta_q^Q(\mathsf C_{\mu,q})
\mid\mu\in\operatorname{OP}_n^G(r)\}$ relative to
$\{\varepsilon_{\lambda,q}^Q
\mid\lambda\in\operatorname{SP}_n(r)\}$ is
$$\bigl(\zeta_q^\lambda(\mu)\bigr)_{\substack{
\mu\in\operatorname{OP}_n^G(r),\\ \lambda\in\operatorname{SP}_n(r)}}\operatorname{diag}_{\lambda\in\operatorname{SP}_n(r)}
\bigl(\kappa_{\lambda,q}\bigr).$$
The matrix $\bigl(\zeta_q^\lambda(\mu)\bigr)_{\substack{
\mu\in\operatorname{OP}_n^G(r),\\ \lambda\in\operatorname{SP}_n(r)}}$ is invertible because its specialization $\mathsf Y_n(r)$ at $v=1$ is invertible by Lemma~\ref{thm:q-mat-class}. Since every $\kappa_{\lambda,q}$ is nonzero, this coefficient matrix is invertible,
and hence the family $\{\eta_q^Q(\mathsf C_{\mu,q})
\mid\mu\in\operatorname{OP}_n^G(r)\}$ is a basis of $Z(\mathcal Q_q(n,r))_{\bar0}$ by
Proposition~\ref{prop:coeff-matrix}.

For the odd family, write
$
\mathsf D_{\lambda,q}
=
\sum_{\rho\in\operatorname{SP}^-(r)}
d_{\lambda\rho}\omega_{\rho,q}.
$
Since $\omega_{\sigma,q}\omega_{\rho,q}
=\delta_{\sigma\rho}e_{\rho,q}^c$, we have
$$
\left(
\frac{\zeta_q^\rho(\mathsf D_{\lambda,q}\mathsf D_{\rho,q})}
     {\zeta_q^\rho(1)}
\right)_{\lambda,\rho\in\operatorname{SP}_n^-(r)}
=
(d_{\lambda\rho})_{\lambda,\rho\in\operatorname{SP}_n^-(r)}
\operatorname{diag}_{\rho\in\operatorname{SP}_n^-(r)}
(d_{\rho\rho}).
$$
By Lemma~\ref{lem:q-casimir-special} and
\cite[Proposition~2.5]{WW13}, the matrix on the left is
regular at $v=1$. Its specialization is an invertible
diagonal matrix by \eqref{eq:Dv=1}.
Hence $(d_{\lambda\rho})_{\lambda,\rho\in\operatorname{SP}_n^-(r)}$
is invertible, and the family
$\{\eta_q^Q(\mathsf D_{\lambda,q})
\mid\lambda\in\operatorname{SP}_n^-(r)\}$
forms a basis of $Z(\mathcal Q_q(n,r))_{\bar1}$ by
Proposition~\ref{prop:coeff-matrix}.
\end{proof}




\end{document}